\documentclass[a4paper,11pt,reqno]{amsart}

\usepackage{amssymb,dsfont}

\chardef\bslash=`\\ % p. 424, TeXbook

\numberwithin{equation}{section}

\newtheorem{theorem}{Theorem}[section]
\newtheorem{corollary}[theorem]{Corollary}
\newtheorem{lemma}[theorem]{Lemma}

\theoremstyle{remark}
\newtheorem{remark}[theorem]{Remark}

\theoremstyle{definition}

\newcommand\bp{\begin{proof}}
\newcommand\ep{\end{proof}}
\newcommand\chf{{\mathds 1}}

\newcommand\aaa{\mathfrak a}
\newcommand\bb{\mathfrak b}
\newcommand\pp{\mathfrak p}
\newcommand\To{\mathfrak T}

\newcommand{\N}{\mathbb N}
\newcommand{\Z}{\mathbb Z}
\newcommand{\Q}{\mathbb Q}
\newcommand{\R}{\mathbb R}

\newcommand\T{\mathbb T}

\newcommand\ak{{\mathbb A}_K}
\newcommand\akf{{\mathbb A}_{K,f}}

\newcommand\OO{{\mathcal O}}
\newcommand\OOh{\widehat{\mathcal O}}
\newcommand\ohs{{\widehat{\OO}^*}}
\newcommand\kab{K^{ab}}
\newcommand\Ind{\operatorname{Ind}}

\newcommand{\Ad}{\operatorname{Ad}}

\newcommand\enu[1]{\smallskip\newline\makebox[6mm][l]{\rm(#1)}}

\begin{document}

\title{KMS states on the C$^*$-algebras of non-principal groupoids}

\author[S. Neshveyev]{Sergey Neshveyev}

\address{Sergey Neshveyev,  Department of Mathematics\\
University of Oslo\\
PO Box 1053 Blindern\\
N-0316 Oslo\\
Norway}
\email{sergeyn@math.uio.no}

\thanks{Partially supported by the Research Council of Norway}

\begin{abstract}
We describe KMS-states on the C$^*$-algebras of etale groupoids in terms of measurable fields of traces on the C$^*$-algebras of the isotropy groups. We use this description to analyze tracial states on the transformation groupoid C$^*$-algebras and to give a short proof of recent results of Cuntz, Deninger and Laca on the Toeplitz algebras of the $ax+b$ semigroups of the rings of integers in number fields.
\end{abstract}

\date{June 29, 2011; minor corrections October 18, 2013.\\ 
This version is different in a few lines from the published version, correcting a potentially misleading omission of the Radon-Nikodym cocycle in the proof of Theorem~\ref{tmain1}; September 23, 2014}

\maketitle

\section*{Introduction}

The problem of classifying KMS-states for various C$^*$-dynamical systems has been extensively studied since the 1970s. Although it can be approached from different angles, one general result is particularly useful and can be applied to almost all known examples. It is the theorem proved  by Renault~\cite{Re1} which states that for the C$^*$-algebra of an etale principal groupoid $G$ with the dynamics given by an $\R$-valued $1$-cocycle $c$, there is a one-to-one correspondence between KMS$_\beta$-states and quasi-invariant probability measures on~$G^{(0)}$ with Radon-Nikodym cocycle $e^{-\beta c}$. Once a C$^*$-algebra is written as a groupoid C$^*$-algebra, this theorem allows one, in many cases, to reduce the study of KMS-states to a measure-theoretic problem, to which results and methods of the dynamical systems theory can be applied; see~e.g.~\cite{KR}. Recently, however, several natural examples of C$^*$-algebras of non-principal groupoids have emerged where the structure of KMS-states is relatively simple, but it cannot be understood in terms of quasi-invariant measures only. One such example is the Toeplitz algebra of the semigroup $\N\rtimes\N^\times$ studied by Laca and Raeburn~\cite{LRt}. In this example one has, for every $\beta>2$, a unique quasi-invariant measure on~$G^{(0)}$ with the required Radon-Nikodym cocycle, but the simplex of KMS$_\beta$-states is isomorphic to the simplex of probability measures on the unit circle. The reason for this structure is that for non-principal groupoids possible extensions of a state on~$C_0(G^{(0)})$ to a KMS-state on~$C^*(G)$ are determined by a choice of tracial states on the C$^*$-algebras of the isotropy groups. In the case studied by Laca and Raeburn these isotropy groups turn out to be~$\Z$, and this explains why measures on the circle appear naturally in the classification of KMS-states. This idea was briefly outlined in the preliminary version~\cite{Ntr} of this note and in the introduction to~\cite{LNt}. In this extended version we formulate the result more explicitly and in a more general setting than discussed in \cite{Ntr,LNt}, and consider more examples.

The proposed strategy for classifying KMS-states on the C$^*$-algebra $C^*(G)$ of an etale groupoid can therefore be described as follows, see Section~\ref{s1} for precise statements. First we have to find all probability measures $\mu$ on~$G^{(0)}$ with Radon-Nikodym cocycle $e^{-\beta c}$. If the $\mu$-measure of the set of points with non-trivial isotropy is zero, then the only way to extend~$\mu_*$ to a KMS$_\beta$-state is by composing~$\mu_*$ with the canonical conditional expectation $C^*(G)\to C_0(G^{(0)})$. Otherwise all possible extensions of~$\mu_*$ are obtained by choosing tracial states $\varphi_x$ on the C$^*$-algebras $C^*(G^x_x)$ of the isotropy groups. The additional requirements on~$\varphi_x$ are that the field $(\varphi_x)_x$ is essentially $G$-invariant and $\mu$-measurable.

These requirements imply that it suffices to specify $\varphi_x$ on a subset intersecting almost every orbit of points with non-trivial isotropy. If the action of $G$ on~$G^{(0)}$ has complicated dynamics and a lot of isotropy, this is hardly a simplification and our description of possible extensions of~$\mu_*$ is probably not very useful. But in many examples the measures $\mu$ are concentrated on a set with only countably many points with non-trivial isotropy. In such cases this description allows us to divide classification of KMS-states into two almost disjoint problems: classification of quasi-invariant measures on~$G^{(0)}$ with a given Radon-Nikodym cocycle and classification of tracial states on~$C^*(G^x_x)$ for countably many points.

This note consists of three sections. Section~\ref{s1} contains our main general results. In Section~\ref{s2} we explain what they mean for transformation groupoids and tracial states. In particular, we will show that classification of tracial states on crossed products by abelian groups reduces completely to a measure-theoretic problem. In Section~\ref{s3} we consider the Toeplitz algebras of the $ax+b$ semigroups of the rings of integers in number fields, recently studied by Cuntz, Deninger and Laca~\cite{CDL}, and recover the classification of KMS-states obtained in~\cite{CDL}. Our approach clarifies why some computations in~\cite{CDL} resemble those used in the study of the Bost-Connes systems of number fields. It also explains why certain representations play a prominent role in the construction of KMS-states in~\cite{CDL}.

\bigskip

\section{KMS states on groupoid C$^*$-algebras} \label{s1}

Let $G$ be a locally compact second countable etale groupoid. We denote the unit space of $G$ by~$G^{(0)}$, and the range and source maps $G\to G^{(0)}$ by $r$ and $s$, respectively. Recall that being etale means that $r$ and $s$ are local homeomorphisms. For $x\in G^{(0)}$ put
$$
G^x=r^{-1}(x),\ \ G_x=s^{-1}(x)\ \ \hbox{and}\ \ G^x_x=G^x\cap G_x.
$$
The C$^*$-algebra $C^*(G)$ of the groupoid $G$ is the C$^*$-enveloping algebra of the $*$-algebra $C_c(G)$ with convolution product
$$
(f_1*f_2)(g)=\sum_{h\in G^{r(g)}}f_1(h)f_2(h^{-1}g)
$$
and involution $f^*(g)=\overline{f(g^{-1})}$.

\smallskip

Let $\mu$ be a probability measure on $G^{(0)}$. Assume that for $\mu$-a.e.~$x\in G^{(0)}$ we are given a state~$\varphi_x$ on~$C^*(G^x_x)$. Denote the generators of $C^*(G^x_x)$ by $u_g$, $g\in G^x_x$. We say that the field of states $\{\varphi_x\}_{x\in G^{(0)}}$ is $\mu$-measurable if for every $f\in C_c(G)$ the function
$$
G^{(0)}\ni x\mapsto\sum_{g\in G^x_x}f(g)\varphi_x(u_g)
$$
is $\mu$-measurable; note that this function is always bounded. We do not distinguish between measurable fields which agree for $\mu$-a.e.~$x$.

\smallskip

Recall that the centralizer $A_\varphi$ of a state $\varphi$ on a C$^*$-algebra $A$ is the set of elements $a\in A$ such that $\varphi(ab)=\varphi(ba)$ for all $b\in A$.

\begin{theorem}\label{tmain1}
There is a one-to-one correspondence between states on $C^*(G)$ with the centralizer containing $C_0(G^{(0)})$ and pairs $(\mu,\{\varphi_x\}_x)$ consisting of a probability measure $\mu$ on $G^{(0)}$ and a $\mu$-measurable field of states $\varphi_x$ on $C^*(G^x_x)$. Namely, the state corresponding to $(\mu,\{\varphi_x\}_x)$ is given~by
$$
\varphi(f)=\int_{G^{(0)}}\sum_{g\in G^x_x}f(g)\varphi_x(u_g)d\mu(x)\ \ \hbox{for}\ \ f\in C_c(G).
$$
\end{theorem}

\bp Assume $\varphi$ is a state on $C^*(G)$ with the centralizer containing $C_0(G^{(0)})$. Let $(H,\pi,\xi)$ be the corresponding GNS-triple. By Renault's disintegration theorem~\cite{Re1,Re2} the representation~$\pi$~is the integrated form of a representation of $G$ on a measurable, with respect to a measure class~$[\nu]$ on~$G^{(0)}$, field of Hilbert spaces $H_x$, $x\in G^{(0)}$. Identifying $H$ with~$\int^\oplus_{G^{(0)}}H_x\,d\nu(x)$, consider the vector field~$(\xi_x)_x$ that defines $\xi$. Let $\mu$ be the measure on $G^{(0)}$ such that $d\mu(x)=\|\xi_x\|^2d\nu(x)$. In other words,~$\mu$~is the measure defined by the restriction of $\varphi$ to $C_0(G^{(0)})$. The action of $G$ on $(H_x)_x$ defines, for every $x$, a representation $\rho_x\colon C^*(G^x_x)\to B(H_x)$. For every $x$ with $\xi_x\ne0$ denote by $\varphi_x$ the state $\|\xi_x\|^{-2}(\rho_x(\cdot)\xi_x,\xi_x)$ on $C^*(G^x_x)$.

For every $f\in C_c(G)$ we have
$$
\varphi(f)=(\pi(f)\xi,\xi)=\int_X\sum_{g\in G^x}D(g)^{-1/2}f(g)(g\xi_{s(g)},\xi_x)d\nu(x),
$$
where $D$ is the Radon-Nikodym cocycle defined by the quasi-invariant measure $\nu$.
Therefore to prove the identity in the formulation of the theorem we have to show that for $\nu$-a.e.~$x$ and every $g\in G^x\setminus G^x_x$ we have $(g\xi_{s(g)},\xi_x)=0$ (note that $D(g)=1$ for $\nu$-a.e.~$x\in G^{(0)}$ and all $g\in G^x_x$). Choose an open subset $U\subset G\setminus G'$, where $G'=\cup_xG^x_x$ is the isotropy bundle, such that $r(U)\cap s(U)=\emptyset$. Let $f\in C_c(U)$. For any function $h\in C_c(r(U))$ we have $f*h=0$. Since $h$ is in the centralizer of $\varphi$, we therefore get
$$
0=\varphi(h*f)=\int_{r(U)}h(x)\sum_{g\in G^x}D(g)^{-1/2}f(g)(g\xi_{s(g)},\xi_x)d\nu(x).
$$
Hence $\sum_{g\in G^x}D(g)^{-1/2}f(g)(g\xi_{s(g)},\xi_x)=0$ for $\nu$-a.e.~$x\in r(U)$. It follows that $(g\xi_{s(g)},\xi_x)=0$ for $\nu$-a.e.~$x\in\nolinebreak r(U)$ and all $g\in G^x\cap U$. Since $G\setminus G'$ can be covered by countably many such sets~$U$, we conclude that  $(g\xi_{s(g)},\xi_x)=0$ for $\nu$-a.e.~$x\in G^{(0)}$ and all $g\in G^x\setminus G^x_x$.

\smallskip

Conversely, assume we are given a probability measure $\mu$ on $G^{(0)}$ and a $\mu$-measurable field of states $\varphi_x$ on $C^*(G^x_x)$. For every $x$ define a state $\psi_x$ on $C^*(G)$ by
$$
\psi_x(f)=\sum_{g\in G^x_x}f(g)\varphi_x(u_g)\ \ \hbox{for}\ \ f\in C_c(G).
$$
In order to show that $\psi_x$ is indeed a well-defined state, consider the GNS-triple $(K_x,\pi_x,\xi_x)$ defined by $\varphi_x$. Induce $\pi_x$ to a representation of $G$ and denote by $\vartheta_x$ its integrated form. Explicitly, $\vartheta_x$ is the representation on the space $L_x$ of functions $\xi\colon G_x\to K_x$ such that
$$
\xi(gh)=\pi_x(u_h^*)\xi(g)\ \ \hbox{for}\ \ g\in G_x\ \ \hbox{and}\ \ h\in G^x_x,\ \ \hbox{and}\ \ \sum_{g\in G_x/G^x_x}\|\xi(g)\|^2<\infty,
$$
given by
$$
(\vartheta_x(f)\xi)(g)=\sum_{h\in G^{r(g)}}f(h)\xi(h^{-1}g)\ \ \hbox{for}\ \ f\in C_c(G).
$$
Let $\zeta_x\in L_x$ be the vector defined by $\zeta_x(g)=\pi_x(u_g^*)\xi_x$ if $g\in G^x_x$ and $\zeta_x(g)=0$ if $\zeta\in G_x\setminus G^x_x$. Then $\psi_x=(\vartheta_x(\cdot)\zeta_x,\zeta_x)$.

Clearly, $C_0(X)$ is contained in the centralizer of $\psi_x$; in fact, for any $f\in C_0(X)$ and $a\in C^*(G)$ we have $\psi_x(fa)=\psi_x(af)=f(x)\psi_x(a)$.
By assumption, the map $x\mapsto\psi_x(a)$ is $\mu$-measurable for every $a\in C_c(G)$, hence for every $a\in\nolinebreak C^*(G)$. Therefore we can define $\varphi(a)=\int_X\psi_x(a)d\mu(x)$.

\smallskip

Finally, it is easy to see that if $(\mu,\{\varphi_x\}_x)$ and $(\tilde\mu,\{\tilde\varphi_x\}_x)$ define the same state then $\mu=\tilde\mu$ and $\varphi_x=\tilde\varphi_x$ for $\mu$-a.e.~$x$.
\ep

Let $E\colon C^*(G)\to C_0(G^{(0)})$ be the canonical conditional expectation. Given a probability measure~$\mu$ on $G^{(0)}$, consider the $\mu$-measurable field of states consisting of canonical traces on $C^*(G^x_x)$. The corresponding state on $C^*(G)$ is $\mu_*\circ E$.

\begin{corollary} \label{ccond}
Let $\mu$ be a probability measure on $G^{(0)}$. Assume that the points of $G^{(0)}$ with non-trivial isotropy form a set of $\mu$-measure zero. Then $\varphi=\mu_*\circ E$ is the unique state on $C^*(G)$ such that the centralizer of $\varphi$ contains $C_0(G^{(0)})$ and $\varphi|_{C_0(G^{(0)})}=\mu_*$.
\end{corollary}

A short proof of this corollary, not relying on the disintegration theorem, can be obtained along the same lines as the proof of~\cite[Proposition~1.1]{LLNlat}.

\smallskip

Let $c$ be a continuous $\R$-valued $1$-cocycle on $G$, that is, a continuous homomorphism $c\colon G\to\nolinebreak\R$. It defines a one-parameter group of automorphisms of $C^*(G)$ by $\sigma^c_t(f)(g)=e^{it c(g)}f(g)$.

Recall that a measure $\mu$ on $G^{(0)}$ is called quasi-invariant with Radon-Nikodym cocycle $e^{c}$ if
$d\mu_r/d\mu_s=e^{c}$, where the measures $\mu_r$ and $\mu_s$ on $G$ are defined by
$$
\int_Gfd\mu_r=\int_{G^{(0)}}\sum_{g\in G^x}f(g)d\mu(x), \ \ \int_Gfd\mu_s=\int_{G^{(0)}}\sum_{g\in G_x}f(g)d\mu(x)\ \ \hbox{for}\ \ f\in C_c(G).
$$
Equivalently, for every open set $U\subset G$ such that $r$ and $s$ are injective on $U$ we have
$$
\frac{d\,T_*\mu}{d\mu}(x)=e^{c(g_x)}\ \ \hbox{for}\ \ x\in s(U),
$$
where $g_x\in U$ is the unique element such that $s(g_x)=x$ and where $T\colon r(U)\to s(U)$ is the homeomorphism defined by $T(r(g_x))=x$.

\smallskip

Now recall that if $\sigma$ is a strongly continuous one-parameter group of automorphisms of a C$^*$-algebra~$A$ and~$\beta\in\R$, then a state $\varphi$ on $A$ is called a $\sigma$-KMS$_\beta$-state if $\varphi$ is $\sigma$-invariant and $\varphi(ab)=\varphi(b\sigma_{i\beta}(a))$ for a dense set of $\sigma$-analytic elements $a,b\in A$.

\begin{theorem} \label{tmain2}
Let $c$ be a continuous $\R$-valued $1$-cocycle on $G$, $\sigma^c$ the dynamics on $C^*(G)$ defined by~$c$, and $\beta\in \R$. Then there exists a one-to-one correspondence between $\sigma^c$-KMS$_\beta$-states on~$C^*(G)$ and pairs $(\mu,\{\varphi_x\}_x)$ consisting of a probability measure $\mu$ on $G^{(0)}$ and a $\mu$-measurable field of states~$\varphi_x$ on $C^*(G^x_x)$ such that:
\enu{i} $\mu$ is quasi-invariant with Radon-Nikodym cocycle $e^{-\beta c}$;
\enu{ii} $\varphi_x(u_g)=\varphi_{r(h)}(u_{hgh^{-1}})$ for $\mu$-a.e.~$x$ and all $g\in G^x_x$ and $h\in G_x$; in particular, $\varphi_x$ is tracial for $\mu$-a.e.~$x$;
\enu{iii} $\varphi_x(u_g)=0$ for $\mu$-a.e.~$x$ and all $g\in G^x_x\setminus c^{-1}(0)$.
\end{theorem}

Note that if $\beta\ne0$ and $\mu$ is a quasi-invariant measure with Radon-Nikodym cocycle $e^{-\beta c}$ then $G^x_x\subset c^{-1}(0)$ for $\mu$-a.e.~$x$. Therefore condition (iii) is redundant in this case, but it is still useful to keep it in mind. On the level of KMS-states this corresponds to the following fact: if $\varphi$ is a state such that $\varphi(ab)=\varphi(b\sigma_{i\beta}(a))$ for $\beta\ne0$, then $\varphi$ is automatically $\sigma$-invariant.

\bp[Proof of Theorem~\ref{tmain2}] Since the centralizer of any $\sigma^c$-KMS-state $\varphi$ contains $C_0(G^{(0)})$, by Theorem~\ref{tmain1} any such state is defined by a pair $(\mu,\{\varphi_x\}_x)$ consisting of a probability measure $\mu$ on $G^{(0)}$ and a $\mu$-measurable field of states~$\varphi_x$ on $C^*(G^x_x)$. Therefore we just have to check that given such a pair $(\mu,\{\varphi_x\}_x)$ the corresponding state is a KMS$_\beta$-state if and only if conditions (i)-(iii) are satisfied.

\smallskip

It is easy to see that condition (iii) is equivalent to $\sigma^c$-invariance.
We will show that (i) and (ii) together are equivalent to $\varphi(f_1*f_2)=\varphi(f_2*\sigma^c_{i\beta}(f_1))$ for all $f_1,f_2\in C_c(G)$.

\smallskip

Assume the equality $\varphi(f_1*f_2)=\varphi(f_2*\sigma^c_{i\beta}(f_1))$ holds for all $f_1,f_2\in C_c(G)$. Fix an open set $U\subset G$ such that $r$ and $s$ are injective on $U$. Let $x\mapsto h^x$ be the inverse of $r\colon U\to r(U)$, and $x\mapsto h_x$ the inverse of $s\colon U\to s(U)$. Denote by $T\colon r(U)\to s(U)$ the homeomorphism defined by $Tx=s(h^x)$, so that $h_{Tx}=h^x$. For $f_1\in C_c(U)$ and $f_2\in C_c(G)$ we have
\begin{align*}
(f_1*f_2)(g)&=\begin{cases}f_1(h^x)f_2((h^x)^{-1}g),&\hbox{if}\ \ x=r(g)\in r(U),\\0,&\hbox{if}\ \ r(g)\notin r(U),\end{cases}\\
(f_2*\sigma^c_{i\beta}(f_1))(g)&=\begin{cases}e^{-\beta c(h_x)}f_1(h_x)f_2(gh_x^{-1}),&\hbox{if}\ \ x=s(g)\in s(U),\\0,&\hbox{if}\ \ s(g)\notin s(U).\end{cases}
\end{align*}
Therefore the equality $\varphi(f_1*f_2)=\varphi(f_2*\sigma^c_{i\beta}(f_1))$ reads as
\begin{equation}\label{eKMSmain}
\int_{r(U)}f_1(h^x)\sum_{g\in G^x_x}f_2((h^x)^{-1}g)\varphi_x(u_g)d\mu(x)
=\int_{s(U)}e^{-\beta c(h_x)}f_1(h_x)\sum_{g\in G^x_x}f_2(gh_x^{-1})\varphi_x(u_g)d\mu(x).
\end{equation}
Let $f\in C_c(s(U))$. Apply the above identity to the functions $f_1$ and $f_2$ defined by $f_1(h_x)=f(x)$ for $x\in s(U)$ and $f_2=f^*_1$. Since $f_1(h^x)=f_1(h_{Tx})=f(Tx)$ for $x\in r(U)$, we get
$$
\int_{r(U)}|f(Tx)|^2d\mu(x)
=\int_{s(U)}e^{-\beta c(h_x)}|f(x)|^2d\mu(x).
$$
Since this is true for all $U$ on which $r$ and $s$ are injective and all $f\in C_c(s(U))$, we see that $\mu$ is quasi-invariant with Radon-Nikodym cocycle $e^{-\beta c}$, so condition (i) is satisfied. But then~\eqref{eKMSmain}, for arbitrary $f_1\in C_c(U)$ and $f_2\in C_c(G)$, can be written as
$$
\int_{r(U)}f_1(h^x)\sum_{g\in G^x_x}f_2((h^x)^{-1}g)\varphi_x(u_g)d\mu(x)
=\int_{r(U)}f_1(h_{Tx})\sum_{g\in G^{Tx}_{Tx}}f_2(gh_{Tx}^{-1})\varphi_{Tx}(u_g)d\mu(x).
$$
Using that $h_{Tx}=h^x$ and $G^x_xh^x=h^xG^{Tx}_{Tx}$, this, in turn, can be written as
$$
\int_{r(U)}f_1(h^x)\sum_{g\in G^x_x}f_2((h^x)^{-1}g)\varphi_x(u_g)d\mu(x)
=\int_{r(U)}f_1(h^x)\sum_{g\in G^{x}_{x}}f_2((h^x)^{-1}g)\varphi_{Tx}(u_{(h^x)^{-1}gh^x})d\mu(x).
$$
Since this equality holds for all $f_1\in C_c(U)$ and $f_2\in C_c(G)$, we conclude that for $\mu$-a.e.~$x\in r(U)$ we have $\varphi_x(u_g)=\varphi_{Tx}(u_{(h^x)^{-1}gh^x})$ for all $g\in G^x_x$. As $G$ can be covered by countably many open sets~$U$ such that $r$ and $s$ are injective on $U$, it follows that $\varphi_x(u_g)=\varphi_{r(h)}(u_{hgh^{-1}})$ for $\mu$-a.e.~$x$ and all $g\in G^x_x$ and $h\in G_x$, so condition (ii) is also satisfied.

\smallskip

Conversely, if conditions (i) and (ii) are satisfied, then we see from the above computations that $\varphi(f_1*f_2)=\varphi(f_2*\sigma^c_{i\beta}(f_1))$ for all $f_1\in C_c(U)$ and $f_2\in C_c(G)$, where $U\subset G$ is any open set such that $r$ and $s$ are injective on $U$. Hence $\varphi(f_1*f_2)=\varphi(f_2*\sigma^c_{i\beta}(f_1))$ for all $f_1,f_2\in C_c(G)$.
\ep

Given a probability measure $\mu$ on $G^{(0)}$ with Radon-Nikodym cocycle $e^{-\beta c}$, the simplest way to extend the state $\mu_*$ to a $\sigma^c$-KMS$_\beta$-state on $C^*(G)$ is by composing $\mu_*$ with the conditional expectation $E\colon C^*(G)\to C_0(G^{(0)})$. As we already noted, this corresponds to taking the canonical traces on~$C^*(G^x_x)$ for~$\varphi_x$. If the points with non-trivial isotropy have measure zero, then $\mu_*\circ E$ is the unique $\sigma^c$-KMS$_\beta$-state extending $\mu_*$. In general, to extend $\mu_*$ we have to find measurable fields of tracial states satisfying properties (ii) and (iii) above. If there are many points with non-trivial isotropy, this is a difficult problem and it is not clear how useful the description in terms of fields of traces is. A simple example of such a rich isotropy bundle structure is the transformation groupoid of the action of the infinite symmetric group $S_\infty$ on $\{0,1\}^\infty$. In many cases, however, the measure~$\mu$~is concentrated on a set with only countably many points with non-trivial isotropy, and then the above result gives a complete description of possible extensions of $\mu_*$.

\begin{corollary} \label{lcount}
Suppose $\mu$ is a quasi-invariant probability measure on $G^{(0)}$ with Radon-Nikodym cocycle~$e^{-\beta c}$. Assume there exists a sequence $\{O_n\}^N_{n=1}$, $N\in\N\cup\{+\infty\}$, of different orbits of the action of $G$ on $G^{(0)}$ such that $\mu(O_n)>0$ for every $n$, and almost all points in $G^{(0)}\setminus\cup_nO_n$ have trivial isotropy. Choose a point $x_n\in O_n$ for every $n$. Then there is a one-to-one correspondence between $\sigma^c$-KMS$_\beta$-states on~$C^*(G)$ extending the state $\mu_*$ on~$C_0(G^{(0)})$ and sequences of tracial states~$\tau_n$ on~$C^*(G^{x_n}_{x_n})$ such that $\tau_n(u_g)=0$ for every $g\in G^{x_n}_{x_n}\setminus c^{-1}(0)$.
\end{corollary}

\bp Every orbit $O_n$ is a countable set, so the measurability assumption is satisfied for any choice of states $\varphi_x$ on $C^*(G^x_x)$ for $x\in O_n$. Since the action of $G$ on $O_n$ is transitive, the map $\{\varphi_x\}_{x\in O_n}\mapsto \varphi_{x_n}$ is a bijection between sequences of states such that $\varphi_x(u_g)=\varphi_{r(h)}(u_{hgh^{-1}})$ for every $x\in O_n$ and all $g\in G^x_x$ and $h\in G_x$ and the set of tracial states on $C^*(G^{x_n}_{x_n})$.
\ep

\bigskip

\section{Traces on crossed products} \label{s2}

Let $X$ be a locally compact second countable space. Assume a countable group $\Gamma$ acts on $X$ by homeomorphisms. Then $C_0(X)\rtimes\Gamma$ is the C$^*$-algebra of the transformation groupoid $X\times\Gamma$. We will apply the results of the previous section to describe tracial states on $C_0(X)\rtimes\Gamma$.

\smallskip

For $x\in X$, denote by $\Gamma_x$ the stabilizer of $x$ in $\Gamma$. Let $\mu$ be a probability measure on $X$. According to our definition, a field of states $\varphi_x$ on $C^*(\Gamma_x)$ is $\mu$-measurable if, for every $g\in\Gamma$, the function $x\mapsto\varphi_x(u_g)$ is $\mu$-measurable on the closed set of points fixed by $g$. This can also be formulated as follows.

Every state $\varphi_x$ extends to a state $\psi_x$ on $C^*(\Gamma)$ such that $\psi_x(u_g)=0$ for every $g\notin\Gamma_x$. This can be proved using induction, exactly as in the proof of Theorem~\ref{tmain1}, so we will write $\Ind^\Gamma_{\Gamma_x}\varphi_x$ for~$\psi_x$. Then $\{\varphi_x\}_x$ is $\mu$-measurable if and only if the map $x\mapsto\Ind^\Gamma_{\Gamma_x}\varphi_x$ from $X$ into the state space $S(C^*(\Gamma))$ of $C^*(\Gamma)$ is $\mu$-measurable, where $S(C^*(\Gamma))$ is considered with the Borel structure defined by the weak$^*$ topology.

Applying Theorem~\ref{tmain2} to the transformation groupoid $X\times\Gamma$ and the zero cocycle, we get the following result.

\begin{theorem} \label{ttrace}
There is a one-to-one correspondence between tracial states on $C_0(X)\rtimes\Gamma$ and pairs $(\mu,\{\varphi_x\}_x)$ consisting of a probability measure $\mu$ on $X$ and a $\mu$-measurable field of states~$\varphi_x$ on~$C^*(\Gamma_x)$ such that:
\enu{i} $\mu$ is $\Gamma$-invariant;
\enu{ii} $\varphi_{x}(u_g)=\varphi_{hx}(u_{hgh^{-1}})$ for $\mu$-a.e.~$x$ and all $g\in \Gamma_x$ and $h\in\Gamma$.
\end{theorem}

Equivalently, we can say that to define a tracial state we need a $\Gamma$-invariant probability measure~$\mu$ on $X$ and a field of tracial states $\varphi_x$ on $C^*(\Gamma_x)$ such that the map $x\mapsto \Ind^\Gamma_{\Gamma_x}\varphi_x\in S(C^*(\Gamma))$ is $\mu$-measurable and $\Gamma$-equivariant, where $g\in\Gamma$ acts on $S(C^*(\Gamma))$ by mapping a state $\psi$ into $\psi\circ(\Ad u_{g})^{-1}$.

\smallskip

The above result was obtained using Renault's disintegration theorem. In the case of transformation groupoids this theorem is a simple consequence of standard results on disintegration of representations of C$^*$-algebras, see~e.g.~\cite[Chapter~IV]{Tak1}. But, in fact, in this case the groupoid picture can be bypassed altogether. In order to show this we will need the following observation, which in a way goes back to \cite[Appendix]{MvN}.

\begin{lemma} \label{lcent}
For any state $\varphi$ on a unital C$^*$-algebra $A$, there exists a unique state $\Phi$ on $A_\varphi^{op}\otimes_{max}A$ such that $\Phi(a\otimes b)=\varphi(ab)$ for all $a\in A_\varphi$ and $b\in A$.
\end{lemma}

\bp We may assume that $A\subset B(H)$ and $\varphi$ is defined by a cyclic vector~$\xi\in H$. Assume first that~$\xi$ is separating for $A''$. Let $J$ be the corresponding modular conjugation. Define a representation $\pi$ of $A_\varphi^{op}\otimes_{max}A$ on $H$ by
$\pi(a\otimes b)=Ja^*Jb=bJa^*J$. If $a\in A_\varphi$, then~$a$ commutes with the modular operator, hence $Ja\xi=a^*\xi$.
Therefore $\pi(a\otimes b)\xi=ba\xi$, so that $\Phi:=(\pi(\cdot),\xi,\xi)$ is the required state.

\smallskip

In the general case we will show that a representation $\pi$ of $A_\varphi^{op}\otimes_{max}A$ such that $\pi(a\otimes b)\xi=ba\xi$ always exists. For $a\in A_\varphi$ and $b,c\in A$ we have
$$
(ba\xi,c\xi)=(b\xi,ca^*\xi).
$$
It follows that for every $a\in A_\varphi$ there exists a well-defined operator $\rho(a)$ on $A\xi$ such that $\rho(a)b\xi=ba\xi$, and then $\rho$ is a representation of $A^{op}_\varphi$ on $A\xi$. Since $\rho(u)$ is unitary for unitary $u$, this is a representation by bounded operators, so it extends to a representation of $A^{op}_\varphi$ on $H$. Its image commutes with $A$, so we can define a representation $\pi$ of $A_\varphi^{op}\otimes_{max}A$ on $H$ by $\pi(a\otimes b)=\rho(a)b$. Then $\pi(a\otimes b)\xi=ba\xi$.
\ep

Assume now that $\varphi$ is a state on $A=C_0(X)\rtimes\Gamma$ with the centralizer containing $C_0(X)$. Using the above lemma define a state $\Phi$ on $C_0(X)\otimes A$. Denote by $j$ the canonical homomorphism $C^*(\Gamma)\to M(A)$. Extending $\Phi$ to the multiplier algebra and composing this extension with $$\operatorname{id}\otimes j\colon C_0(X)\otimes C^*(\Gamma)\to C_0(X)\otimes M(A),$$ we get a state $\Psi$ on $C_0(X)\otimes C^*(\Gamma)$ such that
\begin{equation*}\label{etrace}
\Psi(f\otimes a)=\varphi(fj(a))\ \ \hbox{for}\ \ f\in C_0(X)\ \ \hbox{and}\ \ a\in C^*(\Gamma).
\end{equation*}
Disintegrating $\Psi$ with respect to $C_0(X)$ we get a probability measure $\mu$ on $X$ and states $\psi_x$ on~$C^*(\Gamma)$ such that $\Psi=\int^\oplus_X\psi_x\,d\mu(x)$, that is,
\begin{equation*}
\Psi(f\otimes u_g)=\int_Xf(x)\psi_x(u_g)d\nu(x)\ \ \hbox{for}\ \ f\in C_0(X)\ \ \hbox{and}\ \ g\in \Gamma.
\end{equation*}
Therefore $\varphi_x=\psi_x|_{C^*(\Gamma_x)}$ are exactly the states given by Theorem~\ref{tmain1}.

\smallskip

The above argument suggests yet another way of looking at pairs $(\mu,\{\varphi_x\}_x)$ consisting of a probability measure $\mu$ and a $\mu$-measurable field of states on $C^*(\Gamma_x)$: any such pair defines a state $\Psi=\int^\oplus_X\Ind^{\Gamma}_{\Gamma_x}\varphi_x\,d\mu(x)$ on $C_0(X)\otimes C^*(\Gamma)$, and this way we get all states $\Psi=\int^\oplus_X\psi_x\,d\mu(x)$ such that $\psi_x(u_g)=0$ for $\mu$-a.e.~$x\in X$ and all $g\notin \Gamma_x$.

\smallskip

For abelian groups this point of view combined with Theorem~\ref{ttrace} allow us to completely reduce the classification of tracial states to a measure-theoretic problem.

\begin{corollary} \label{cabel}
If $\Gamma$ is abelian, there is a bijection between tracial states $\tau$ on $C_0(X)\rtimes \Gamma$ and probability measures $\nu=\int^\oplus_X\nu_x\,d\mu(x)$ on~$X\times\hat \Gamma$ such that
\enu{i} $\nu$ is invariant with respect to the action of $\Gamma$ on the first factor of $X\times\hat \Gamma$; equivalently, $\mu$ is $\Gamma$-invariant and $\nu_{x}=\nu_{gx}$ for $\mu$-a.e.~$x\in X$ and every $g\in\Gamma$;
\enu{ii} $\nu_x$ is $\Gamma_x^\perp$-invariant for $\mu$-a.e.~$x\in X$.

\smallskip\noindent
Namely, the trace $\tau$ corresponding to such a measure $\nu$ is given by
$$
\tau(fu_g)=\int_{X\times\hat \Gamma}f(x)\chi(g)d\nu(x,\chi) \ \ \hbox{for}\ \ f\in C_0(X)\ \ \hbox{and}\ \ g\in \Gamma.
$$
\end{corollary}

\bp We only need to note that if $\varphi_x$ is the state on $C^*(\Gamma)=C(\hat\Gamma)$ defined by $\nu_x$ then $\varphi_x(u_g)=0$ for all $g\notin\Gamma_x$ if and only if $\nu_x$ is $\Gamma_x^\perp$-invariant.
\ep

Note that in this case it is particularly easy to check that the state $\tau$ corresponding to $\nu$ exists. Indeed, define a representation~$\rho$ of $C_0(X)\rtimes\Gamma$ on~$L^2(X\times\hat \Gamma,d\nu)$ by
$$
(\rho(f)\zeta)(x,\chi)=f(x)\zeta(x,\chi),\ \ (\rho(u_g)\zeta)(x,\chi)=\chi(g)\zeta(g^{-1}x,\chi)
$$
and consider the function $\xi\equiv1$ in $L^2(X\times\hat G,d\nu)$. Then $\tau=(\rho(\cdot)\xi,\xi)$.

\begin{corollary}
If $\Gamma$ is abelian, there exists a one-to-one correspondence between extremal tracial states on $C_0(X)\rtimes\Gamma$ and triples $(H,\chi,\mu)$, where $H$ is a subgroup of $\Gamma$, $\chi$ is a character of $H$ and~$\mu$~is an ergodic $\Gamma$-invariant probability measure $\mu$ on~$X$ such that $\Gamma_x=H$ for $\mu$-a.e.~$x\in X$. Namely, the trace corresponding to $(H,\chi,\mu)$ is given by
$$
\tau(fu_g)=\begin{cases}\chi(g)\int_Xf(x)d\mu(x), &\hbox{if}\ \ g\in H,\\
0,&\hbox{otherwise.}\end{cases}
$$
\end{corollary}

\bp Extremal tracial states correspond to extremal probability measures $\nu$ on $X\times\hat \Gamma$ with properties~(i) and~(ii) in Corollary~\ref{cabel}. If $\nu$ is extremal, then its projection $\mu$ onto $X$ is an extremal $\Gamma$-invariant measure, that is, $\mu$ is ergodic. But then $\nu$ must be of the form $\mu\times\lambda$ for a probability measure $\lambda$ on $\hat\Gamma$. Furthermore, for every $g\in \Gamma$ the set $X_g=\{x\mid gx=x\}$ is $\Gamma$-invariant, so either $X_g$ or its complement has measure zero. In other words, on a subset of full measure we have $\Gamma_x=H$ for a subgroup $H\subset \Gamma$. Hence $\lambda$ is $H^\perp$-invariant, so it can be considered as a measure on $\hat\Gamma/H^\perp=\hat H$. As $\nu$ is extremal, $\lambda$ is an extremal probability measure on $\hat H$, that is, $\lambda=\delta_\chi$ for some $\chi\in\hat H$.

Conversely, any measure $\nu$ of the form $\mu\times\delta_\chi$, where $\mu$ is ergodic and $\chi\in \hat H$, where $H$ is the common stabilizer of the points in a subset of $X$ of full measure, is an extremal measure with properties~(i) and~(ii).
\ep

\begin{remark}
If $\tau$ is the extremal tracial state corresponding to a triple $(H,\chi,\mu)$, then there is a non-canonical isomorphism $\pi_\tau(C_0(X)\rtimes\Gamma)''\cong L^\infty(X,\mu)\rtimes \Gamma/H$. Namely, extend $\chi$ to a character $\tilde\chi$ of $\Gamma$. Then the required isomorphism $\rho$ is given by $\rho(f)=f\in L^\infty(X,\mu)$, $\rho(u_g)=\tilde\chi(g)u_{\bar g}$, where $\bar g$ is the image of $g$ in $\Gamma/H$. Note that since $L^\infty(X,\mu)\rtimes \Gamma/H$ is a factor, this provides a direct proof of the extremality of $\tau$.
\end{remark}

Finally, consider the case $\Gamma=\Z$ explicitly used in~\cite{LNt}. Denote by $T\colon X\to X$ the homeomorphism corresponding to~$1\in\Z$. For $n\ge0$ denote by $X_n\subset X$ the subset of points of period $n$ (so $X_0$ is the set of aperiodic points). Then any measure $\nu$ on $X\times\T$ with properties (i) and (ii) from Corollary~\ref{cabel} decomposes into a sum of measures satisfying the same properties and concentrated on $X_n\times\T$ for some $n$.

If $\nu$ is concentrated on $X_0$ then $\nu=\mu\times\lambda$, where $\lambda$ is the Lebesgue measure, and the corresponding trace is $\mu_*\circ E$, where $E\colon C_0(X)\rtimes \Z\to C_0(X)$ is the canonical conditional expectation.

If $\mu$ is concentrated on $X_n$, $n\ge1$, then $\nu$ is a $\Z/n\Z\times\Z/n\Z$-invariant measure on $X_n\times\T$, where the second factor $\Z/n\Z$ acts on $\T$ by rotations. Consider the simplest case where $\mu$ is concentrated on the orbit of a point~$x$ of period $n$. Then $\mu=n^{-1}\sum^{n-1}_{k=0}\delta_{T^kx}$ and $\nu=\mu\times\lambda$, where $\lambda$ is a measure that is invariant under the rotation by~$2\pi/n$ degrees (we will say that $\lambda$ is $n$-rotation invariant). The corresponding trace can be written as follows.

The $*$-homomorphism $\rho\colon C_0(X)\to C(\Z/n\Z)$, $\rho(f)(k)=f(T^kx)$, extends to a $*$-homomorphism $$\rho\colon C_0(X)\rtimes\Z\to C(\Z/n\Z)\rtimes\Z.$$ Passing to the dual groups we can identify $C(\Z/n\Z)\rtimes\Z$ with $C(\T)\rtimes\Z/n\Z$.  By composing the canonical conditional expectation $C(\T)\rtimes\Z/n\Z\to C(\T)=C^*(\Z)$ with $\rho$ we then get a completely positive map
$$
E_x\colon C_0(X)\rtimes\Z\to C^*(\Z),\ \ E_x(fu^m)=\frac{1}{n}\left(\sum^{n-1}_{k=0}f(T^kx)\right)u^m,
$$
where $u=u_1\in M(C_0(X)\rtimes\Z)$ is the canonical unitary implementing $T$, so $ufu^*=f(T^{-1}\cdot)$. The measure~$\lambda$ defines a state $\lambda_*$ on $C^*(\Z)$ by
$
\lambda_*(u^m)=\int_\T z^md\lambda(z).
$
Then $\lambda_*\circ E_x$ is the required tracial state on~$C_0(X)\rtimes\Z$.

It follows that in the case when there are only countably many periodic orbits Corollary~\ref{cabel} for~$\Gamma=\Z$ can be formulated as follows.

\begin{corollary} \label{ctrace}
Assume a homeomorphism $T$ of $X$ has at most countably many periodic orbits $O_n$. For every $n$ choose $x_n\in O_n$. Then any tracial state $\tau$ on $C_0(X)\rtimes\Z$ has a unique decomposition
$$
\tau=\mu_*\circ E+\sum_{n}{\lambda_n}_*\circ E_{x_n},
$$
where $\mu$ is a $T$-invariant measure on $X$ such that $\mu(O_n)=0$ for every $n$, $\lambda_n$ is an $|O_n|$-rotation invariant measure on $\T$, and $\mu(X)+\sum_n\lambda_n(\T)=1$. Conversely, any such collection of measures~$\mu$, \nolinebreak $\lambda_n$ defines a tracial state.
\end{corollary}

\bigskip

\section{KMS states on the Toeplitz algebras of $ax+b$ semigroups} \label{s3}

In a recent paper Cuntz, Deninger and Laca~\cite{CDL} have studied the Toeplitz algebra of the $ax+b$ semigroup of the ring of integers in a number field. In this section we will show how to recover their classification of KMS-states using our general framework. %As we will see, the analysis is very similar to that in~\cite{LNt} for the Toeplitz algebra of the semigroup $\Z\rtimes\N$, but there are some interesting complications due to the presence of non-trivial units and non-principal ideals in the ring.
In addition to illustrating the general theory, our goal is to clarify a connection between the rather involved analysis in~\cite{CDL} and that of the Bost-Connes systems of number fields.

\smallskip

We will follow the notation in~\cite{LNT} rather than the one in~\cite{CDL}. Let $K$ be a number field with the ring of integers $\OO$. Denote by $V_K$ the set of places of $K$, and by $V_{K,f}\subset V_K$ the subset of finite places. For every place $v$ denote by~$K_v$ the completion of~$K$ at $v$.  For $v\in V_{K,f}$, let $\OO_v$ be the closure of~$\OO$ in~$K_v$. We write $v|\infty$ when~$v$ is infinite and denote by $K_\infty=\prod_{v|\infty}K_v$ the completion of~$K$ at all infinite places. The adele ring $\ak$ is the restricted product, as $v$ ranges over all places, of the rings~$K_v$, with respect to $\OO_v\subset K_v$ for $v\in V_{K,f}$. When the product is taken only over finite places~$v$, we get the ring $\akf$ of finite adeles; we then have $\ak=K_\infty\times\akf$.   The ring of integral adeles is $\OOh=\prod_{v\in V_{K,f}}\OO_v\subset \akf$. Denote by $N_K\colon\akf^*\to(0,+\infty)$ the absolute norm.

\smallskip

Let $\To[\OO]$ be the Toeplitz algebra of the semigroup $\OO\rtimes\OO^\times$~\cite{CDL}, where $\OO^\times$ is the semigroup of nonzero elements in $\OO$. Although $\To[\OO]$ can be defined in terms of generators and relations, we will use a presentation of $\To[\OO]$ as a groupoid C$^*$-algebra. For this consider the space $\Omega_K$ defined as the quotient of $\akf\times\akf/\ohs$ by the equivalence relation
$$
(r,a)\sim (s,b)\ \ \Leftrightarrow\ \  a=b\ \ \hbox{and}\ \ r-s\in a\OOh.
$$
In other words, $\Omega_K$ consists of pairs $(r,a)$ with $a\in \akf/\ohs$ and $r\in\akf/a\OOh$. It is a locally compact space with the quotient topology. Denote by $\Omega_\OO$ the compact open subset of $\Omega_K$ consisting of pairs $(r,a)$ with $a\in \OOh/\ohs$ and $r\in\OOh/a\OOh$. The group $K\rtimes K^*$ acts on $\akf\times\akf$ by $(n,k)(r,a)=(n+kr,ka)$. This action defines an action of $K\rtimes K^*$ on $\Omega_K$. By \cite[Propositions~5.1 and~5.2]{CDL} there is a canonical isomorphism
$$
\To[\OO]\cong \chf_{\Omega_\OO}(C_0(\Omega_K)\rtimes K\rtimes K^*)\chf_{\Omega_\OO}.
$$
Therefore $\To[\OO]$ is the C$^*$-algebra of the reduction of the transformation groupoid $\Omega_K\times(K\rtimes K^*)$ by $\Omega_\OO$. The homomorphism $K\rtimes K^*\to\R$, $(n,k)\mapsto -\log N_K(k)$, defines a $1$-cocycle on the groupoid, which in turn defines a dynamics $\sigma$ on $\To[\OO]$.

According to our general scheme, in order to classify $\sigma$-KMS$_\beta$-states on $\To[\OO]$ we first have to find probability measures $\mu$ on $\Omega_\OO$ with Radon-Nikodym cocycle $N_K^{\beta}$. In other words, we are looking for measures $\mu$ such that
$$
\mu((n,k)Y)=N_K(k)^{\beta}\mu(Y)\ \hbox{for Borel}\ Y\subset\Omega_\OO\ \hbox{and}\ (n,k)\in K\rtimes K^*\ \hbox{such that}\ (n,k)Y\subset\Omega_\OO.
$$
Since any translationally $\OO$-invariant measure on $\OOh$ is a Haar measure, it is not difficult to show, see the proof of \cite[Proposition 2.1]{LNt}, that any measure $\mu$ as above must be the image under the projection $\OOh\times\OOh/\ohs\to \Omega_\OO$ of a measure $m\times\nu$, where $m$ is the normalized Haar measure on~$\OOh$ and~$\nu$~is a probability measure on $\OOh/\ohs$. Furthermore, as $m(k\,\cdot)=N_K(k)m$, the above condition on~$\mu$ is satisfied if and only if
\begin{equation} \label{et1}
\nu(kY)=N_K(k)^{\beta-1}\nu(Y)\ \hbox{for Borel}\ Y\subset\OOh/\ohs\ \hbox{and}\ k\in K^*\ \hbox{such that}\ kY\subset\OOh/\ohs.
\end{equation}
Clearly, there are no such measures for $\beta<1$, so there are no $\sigma$-KMS$_\beta$-states on $\To[\OO]$ for $\beta<1$. The case $\beta=1$ is also easy: the only such measure $\nu$ is concentrated at $0$. We will denote this measure by $\nu_1$. Assume now that $\beta>1$. A more general problem is studied in \cite{LNT}. Namely, we have the following equivalent form of \cite[Theorem 2.5]{LNT}.

\begin{theorem} \label{tlnt}
Let $K^*_+\subset K^*$ be the subgroup of totally positive elements. Consider probability measures $\nu$ on $\OOh$ such that
\begin{equation} \label{et10}
\nu(kY)=N_K(k)^{\beta-1}\nu(Y)\ \hbox{for Borel}\ Y\subset\OOh\ \hbox{and}\ k\in K^*_+\ \hbox{such that}\ kY\subset\OOh.
\end{equation}
Then
\enu{i} for every $1<\beta\le2$ there exists a unique probability measure on $\OOh$ satisfying \eqref{et10};
\enu{ii} for every $\beta>2$ and $a\in\akf^*$ there exists a unique probability measure on $\OOh$ satisfying \eqref{et10} that is concentrated on $\overline{K^*_+}a\cap\OOh$; this way we get, for every $\beta>2$, a one-to-one correspondence between extremal probability measures satisfying \eqref{et10} and points of the compact group $\akf^*/\overline{K^*_+}$; in particular, the simplex of measures satisfying \eqref{et10} is canonically isomorphic to the simplex of probability measures on $\akf^*/\overline{K^*_+}$, and any such measure is concentrated on $\akf^*\cap\OOh$.
\end{theorem}

Here $\overline{K^*_+}$ denotes the closure of $K^*_+$ in $\akf^*$.

We now apply this theorem to classify measures satisfying \eqref{et1}. Clearly, there is a one-to-one correspondence between such measures and $\ohs$-invariant measures satisfying \eqref{et10}.

\smallskip

Let us first consider the case $\beta>2$. We will need only one conclusion from Theorem~\ref{tlnt}(ii): any measure $\nu$ satisfying \eqref{et1} is concentrated on $(\akf^*\cap\OOh)/\ohs$. Recall that $\akf^*/\ohs$ can be identified with the group $J_K$ of fractional ideals, while $K^*/\OO^*$ is the group of principal fractional ideals. Therefore if $h_K$ is the class number of $K$, then the action of $K^*$ by multiplication on~$\akf^*/\ohs$ has exactly $h_K$ orbits. If $\aaa\in\akf^*/\ohs=J_K$, then obviously there exists at most one probability measure satisfying~\eqref{et1} that is concentrated on $K^*\aaa\cap \OOh/\ohs$.
Such a measure indeed exists:
$$
\nu_{\aaa,\beta}=\frac{1}{\zeta(\beta-1,[\aaa])}\sum_{\bb\in J_K^+\cap[\aaa]}N_K(\bb)^{-(\beta-1)}\delta_{\bb},
$$
where $[\aaa]$ is the class of $\aaa$ in the ideal class group $Cl(K)$, $J_K^+\subset J_K$ is the subsemigroup of integral ideals, and $\zeta(\cdot,[\aaa])$ is the partial $\zeta$-function defined by $\zeta(s,[\aaa])=\sum_{\bb\in J_K^+\cap[\aaa]}N_K(\bb)^{-s}$ for $s>1$. Since any measure satisfying \eqref{et1} is concentrated on $J_K^+=(\akf^*\cap\OOh)/\ohs$, we conclude that if $\aaa_1,\dots,\aaa_{h_K}$ are representatives of different ideal classes, then any probability measure satisfying \eqref{et1} is a unique convex combination of the measures $\nu_{\aaa_n,\beta}$, $1\le n\le h_K$.

Let $\aaa\in J_K^+$ be a nonzero ideal in $\OO$. Denote by $\mu_{\aaa,\beta}$ the image of the measure $m\times\nu_{\aaa,\beta}$ under the projection $\OOh\times\OOh/\ohs\to\Omega_\OO$. By construction the measure $\mu_{\aaa,\beta}$ is concentrated on the set
$$
\{(r,\bb)\mid \bb\in J_K^+\cap [\aaa],\ r\in\OOh/\bb\OOh=\OO/\bb\}\subset\Omega_\OO.
$$
The partially defined action of $K\rtimes K^*$ on $\Omega_\OO$ is transitive on this set. Therefore by Corollary~\ref{lcount}, in order to extend $(\mu_{\aaa,\beta})_*$ to a KMS$_\beta$-state we have to choose a tracial state on the C$^*$-algebra of the stabilizer of one point in this set. We take $(0,\aaa)$ as such point. Its stabilizer in $K\rtimes K^*$ is $\aaa\rtimes\OO^*$. For a tracial state $\tau$ on $C^*(\aaa\rtimes\OO^*)$ denote by $\varphi_{\aaa,\tau,\beta}$ the corresponding $\sigma$-KMS$_\beta$-state on $\To[\OO]$.

Now pick representatives $\aaa_1,\dots,\aaa_{h_K}$ of different ideal classes. Any measure $\mu$ on $\Omega_\OO$ with Radon-Nikodym cocycle $N_K^\beta$ is a convex combination of the measures $\mu_{\aaa_n,\beta}$. It is concentrated on the set $\{(r,\aaa)\mid\aaa\in J_K^+,\ r\in\OO/\aaa\}\subset\Omega_\OO$. The partially defined action of $K\rtimes K^*$ has $h_K$ orbits on this set, and we can take the points $(0,\aaa_n)$ as representatives of these orbits. To extend $\mu_*$ we need to choose a tracial state on the C$^*$-algebra of the stabilizer of every point $(0,\aaa_n)$ that carries a positive measure. In other words, any $\sigma$-KMS$_\beta$-state on $\To[\OO]$ is a convex combination $\sum^{h_K}_{n=1}\lambda_n\varphi_{\aaa_n,\tau_n,\beta}$, where $\tau_n$ is a tracial state on the C$^*$-algebra~$C^*(\aaa_n\rtimes\OO^*)$. The weights $\lambda_n$ and the tracial states $\tau_n$ define a tracial state $\oplus_n\lambda_n\tau_n$ on $\oplus_nC^*(\aaa_n\rtimes\OO^*)$. We therefore get the following result.

\begin{theorem}{\cite[Theorem 7.3]{CDL}}
Choose representatives $\aaa_1,\dots,\aaa_{h_K}\in J_K^+$ of different ideal classes. Then, for every $\beta>2$, there is an affine isomorphism between the simplex of $\sigma$-KMS$_\beta$-states on~$\To[\OO]$ and the simplex of tracial states on the C$^*$-algebra $\bigoplus^{h_K}_{n=1}C^*(\aaa_n\rtimes\OO^*)$.
\end{theorem}

Note that the extremal KMS$_\beta$-states are states of the form $\varphi_{\aaa,\tau,\beta}$, where $\tau$ is an extremal tracial state on $C^*(\aaa\rtimes\OO^*)$. Such a state has type I$_\infty$ or II$_\infty$ depending on whether $\tau$ is of type I or II$_1$.

\smallskip

Let us now turn to the more complicated case $1\le\beta\le2$. Recall that the unique measure~$\nu_1$ satisfying \eqref{et1} for $\beta=1$ is the delta-measure at $0\in\OOh/\ohs$. Assume $1<\beta\le2$. Again it is easy to construct a probability measure satisfying~\eqref{et1}: take $\nu_\beta=\prod_{v\in V_{K,f}}\nu_{\beta,v}$, where the measure~$\nu_{\beta,v}$ on~$\OO_v/\OO_v^*$ is defined by
$$
\nu_{\beta,v}=(1-N_K(\pp_v)^{-(\beta-1)})\sum_{n=0}^\infty N_K(\pp_v)^{-(\beta-1)n}\delta_{\pp^n_v},
$$
where $\pp_v$ is the prime ideal in $\OO$ corresponding to the place $v$ and where we identified $\OO^\times_v/\OO_v^*$ with the sequence $\{\pp_v^n\}_n$.\footnote{The measure $\nu_\beta$ is of course well-defined for all $\beta>1$. For $\beta>2$ we have $\displaystyle \nu_\beta=\sum^{h_K}_{n=1}\frac{\zeta(\beta-1,[\aaa_n])}{\zeta_K(\beta-1)}\nu_{\aaa_n,\beta}$, where $\zeta_K$ is the Dedekind $\zeta$-function.}
By Theorem~\ref{tlnt}(i) this is the unique probability measure satisfying \eqref{et1}. Note that the $\nu_\beta$-measure of the set $(\akf^*\cap\OOh)/\ohs$ is zero, which follows from the divergence to zero of the product $\prod_v(1-N_K(\pp_v)^{-(\beta-1)})$. It is also clear that the $\nu_\beta$-measure of the set of points $a\in\OOh/\ohs$ with at least one zero coordinate is zero.

For every $1\le\beta\le 2$, let $\mu_\beta$ be the image of the measure $m\times\nu_\beta$ under the map $\OOh\times\OOh/\ohs\to\Omega_\OO$.

\begin{lemma}
For $1\le\beta\le2$ the set of points in $\Omega_\OO$ with non-trivial stabilizers in $K\rtimes K^*$ has $\mu_\beta$-measure zero.
\end{lemma}

\bp For $\beta=1$ the measure space $(\Omega_\OO,\nu_1)$ can be identified with $(\OOh,m)$, with the partially defined action of $K\rtimes K^*$ given by $(n,k)r=n+kr$. Clearly, every element $g\ne e$ in $K\rtimes K^*$ has at most one fixed point in $\OOh$, and the measure of every such point is zero.

Assume now that $1<\beta\le2$. We have to show that for every element $g=(n,k)\in K\rtimes K^*$, $g\ne e$, the set of points $(r,a)\in\Omega_\OO$ fixed by $g$ has $\mu_\beta$-measure zero. In other words, for $\nu_\beta$-a.e.~$a\in\OOh/\ohs$ the set $A_{g,a}$ of points $r\in\OOh/a\OOh$ such that $g(r,a)=(r,a)$ has $m_a$-measure zero, where $m_a$ is the normalized Haar measure on $\OOh/a\OOh$. We will show that $m_a(A_{g,a})=0$ for every $a=(a_v)_v$ such that $a\notin\akf^*/\ohs$ and $a_v\ne0$ for all $v$. Since the $\nu_\beta$-measure of the complement of the set of such points~$a$ is zero, this will prove the lemma.

The set $A_{g,a}$ can be nonempty only when $ka=a$, whence $k\in\OO^*$ as $a_v\ne0$ for all $v$ by assumption. Then $g$ acts on $\akf/a\OOh$, and the set $A_{g,a}$ consists of points $r\in\OOh/a\OOh$ such that $(k-1)r=-n$. If $k=1$, this means that $A_{g,a}$ is nonempty only when $n\in a\OOh$. But then $g=(0,1)$, since $K\cap a\OOh=\{0\}$ by our assumption that $a\notin\akf^*$. Since we assumed that $g\ne e$, we conclude that $A_{g,a}$ can be nonempty only when $k\ne1$. In this case $k-1$ is invertible in $\OO_v$ for $v$ outside a finite set $F\subset V_{K,f}$. Hence, if $r=(r_v)_v\in A_{g,a}$, then $r_v\in\OO_v/a_v\OO_v$ is uniquely determined for $v\notin F$. It follows that
$$
m_a(A_{g,a})\le\prod_{v\in V_{K,f}\setminus F}|\OO_v/a_v\OO_v|^{-1}.
$$
The latter product diverges to zero, since by assumption there are infinitely many places $v$ such that $a_v\ne\OO^*_v$.
\ep

By Theorem~\ref{tmain2} and Corollary~\ref{ccond} the only way to extend the state $\mu_{\beta*}$ on $C(\Omega_\OO)$ to a KMS$_\beta$-state on $\To[\OO]$ is by composing it with the canonical conditional expectation $E\colon\To[\OO]\to C(\Omega_\OO)$. Thus we have proved the following result.

\begin{theorem}{\cite[Theorem~6.7]{CDL}} \label{tcdl}
For every $1\le\beta\le2$ there exists a unique $\sigma$-KMS$_\beta$-state on $\To[\OO]$.
\end{theorem}

Note that using \cite[Corollary~3.2]{nes2} and the same arguments as in the proof of~\cite[Theorem~3.2]{LNt}, it is easy to show that these KMS-states have type III$_1$.

\smallskip

We would like to finish by making a few remarks about the uniqueness of measures satisfying \eqref{et1} for $1<\beta\le2$, which is the most non-trivial part in the above analysis. When the field $K$ has class number one, this uniqueness is quite simple and a proof can be obtained using the same arguments as for $K=\Q$. This goes back to~\cite{bos-con} and is equivalent to the uniqueness of KMS$_{\beta-1}$-states on the symmetric part of the Bost-Connes system for $\Q$. The situation changes drastically when $h_K>1$. In this case, as we saw, the uniqueness can be deduced from Theorem~\ref{tlnt}(i). Although an equivalent form of this theorem is explicitly stated in \cite{LNT}, it is based on a result on the Bost-Connes systems established in \cite{LLNlat}. Let us briefly describe the arguments in \cite{LNT} applied to the classification problem of measures satisfying \eqref{et1}.

The action of $K^*$ on $\akf/\ohs$ defines an action of the group of principal ideals. Induce this action to an action of the whole group $J_K$ of fractional ideals. By a general result on Morita equivalent systems \cite[Theorem~3.2]{LN}, the set of KMS-weights remains the same under induction. In the present case, when everything is formulated in terms of measures, this is quite obvious: a Radon measure with a given Radon-Nikodym cocycle with respect to a group action is completely determined by its restriction to any compact open subset intersecting every orbit. From this it can be deduced that classification of measures satisfying \eqref{et1} is equivalent to classifying probability measures $\nu$ on $Cl(K)\times \OOh/\ohs$ such that
\begin{equation} \label{et3}
\nu(\aaa Y)=N_K(\aaa)^{-(\beta-1)}\nu(Y)\ \ \hbox{for Borel}\ \ Y\subset Cl(K)\times\OOh/\ohs\ \ \hbox{and}\ \ \aaa\in J_K.
\end{equation}
Here the action of $J_K=\akf^*/\ohs$ on $Cl(K)\times \akf/\ohs$ is diagonal, so for every $v$ the element $\pp_v\in J_K$ changes only two coordinates, one corresponds to $Cl(K)$, the other to the place $v$.

The Bost-Connes system for $K$, in turn, is defined using the partially defined diagonal action of~$J_K$ on the space $$Y_K=Gal(\kab/K)\times_\ohs\OOh.$$ Here $\ohs$ acts on $Gal(\kab/K)$ via the Artin map $r_K\colon\ak^*\to Gal(\kab/K)$. We have $$Gal(\kab/K)/r_K(K^*_\infty\ohs)=Gal(H(K)/K)\cong Cl(K),$$
where $H(K)$ is the Hilbert class field of $K$. Thus
$$
Cl(K)\times \OOh/\ohs=Y_K/r_K(K^*_\infty\ohs).
$$
Therefore the uniqueness of measures satisfying \eqref{et3}, or \eqref{et1}, is equivalent to the uniqueness of $r_K(K^*_\infty\ohs)$-invariant measures on $Y_K$ satisfying the same scaling property. In other words, Theorem~\ref{tcdl} is essentially equivalent to the uniqueness of $r_K(K^*_\infty\ohs)$-invariant KMS$_{\beta-1}$-states on the Bost-Connes system for $K$ for $1<\beta\le2$. That a KMS$_{\beta-1}$-state for the Bost-Connes system is unique for every inverse temperature in this region is proved in \cite[Theorem 2.1]{LLNlat}.

\end{document}